# MÁS ALLÁ DE LA FILA POR COLUMNA: INTERVENCIÓN CON GEOGEBRA Y ETM EN EL APRENDIZAJE DEL PRODUCTO MATRICIAL EN SECUNDARIA

Felix **De la Cruz** Serrano
Institución Educativa María Reiche
Perú
feldese@gmail.com

**Resumen**

La enseñanza del producto de matrices en la educación secundaria suele centrarse en la aplicación mecánica de la regla algorítmica fila × columna, dejando en segundo plano su interpretación como transformación geométrica. Este estudio analiza el impacto de una secuencia didáctica mediada por GeoGebra y fundamentada en el Espacio de Trabajo Matemático (ETM) en el aprendizaje del producto matricial. Participaron diez estudiantes de quinto de secundaria de una institución educativa de Lima (Perú). La intervención se desarrolló en cuatro sesiones, combinando actividades manuales y exploración digital con GeoGebra.

Los resultados evidencian avances significativos en la génesis semiótica, reflejados en la coordinación entre los registros algebraico, gráfico y numérico; en la génesis instrumental, a través del uso progresivamente significativo de GeoGebra como herramienta cognitiva; y en la génesis discursiva, mediante la formulación de explicaciones sobre el efecto geométrico de las matrices. Se observa una transición desde una ejecución algorítmica del producto matricial hacia una comprensión conceptual basada en transformaciones lineales.

Estos hallazgos sugieren que el diseño de tareas que articulan trabajo manual, visualización dinámica y argumentación matemática favorece una comprensión más profunda del producto de matrices y aporta criterios para el uso reflexivo de tecnologías digitales en la enseñanza del álgebra lineal escolar.

**Palabras clave:** Producto de matrices, GeoGebra, ETM, Educación secundaria.

## Introducción

En el currículo de secundaria, la enseñanza del producto de matrices suele reducirse al procedimiento "fila × columna", sin enfatizar su interpretación geométrica como transformación lineal del plano. Este abordaje algorítmico empobrece la comprensión conceptual: el alumnado no consigue vincular cada entrada de la matriz con el estiramiento, la rotación o la reflexión que provoca sobre una figura, ni entender por qué el producto matricial es no conmutativo. Observaciones en aulas de quinto de secundaria confirman que, aunque los estudiantes resuelven operaciones, no explican el significado de las transformaciones ni el papel de cada coeficiente; esta brecha conceptual persiste en niveles posteriores (Hidayanti, 2020).

Frente a tal escenario, GeoGebra ofrece una vía para enriquecer la enseñanza del álgebra matricial, pues permite visualizar en tiempo real el efecto de una matriz sobre un objeto (Turgut, 2022). Sin embargo, la investigación advierte que la mera incorporación de tecnología no basta: sin tareas que obliguen a coordinar los registros algebraico, gráfico y numérico y a fundamentar los resultados, GeoGebra se convierte en un verificador mecánico (Drijvers, 2013).

Para abordar este desafío, el Espacio de Trabajo Matemático (ETM) de Kuzniak (2011) aporta un marco idóneo para analizar la articulación entre registros de representación, herramientas y discurso matemático. Observaciones preliminares en aulas de quinto de secundaria revelaron un ETM fragmentado, caracterizado por la ausencia de herramientas tecnológicas en las prácticas cotidianas y por una escasa articulación entre el cálculo matricial y su interpretación geométrica.

En este contexto, el objetivo general de este estudio es analizar la evolución del Espacio de Trabajo Matemático en estudiantes de secundaria durante el aprendizaje del producto de matrices mediante una secuencia didáctica mediada por GeoGebra. Este manuscrito corresponde a un preprint; una versión preliminar de este trabajo fue presentada previamente como ponencia en un evento académico.

## Elementos teóricos y metodológicos

### Marco teórico

El estudio se apoya en el Espacio de Trabajo Matemático (ETM), modelo que describe cómo se organiza la actividad matemática escolar y que ha sido formalizado por Kuzniak (2011). El ETM distingue dos planos complementarios:
- Epistemológico: reúne los saberes, los registros de representación y los artefactos disponibles para abordar un problema;
- Cognitivo: engloba los procesos de visualización, construcción y validación que desarrolla el estudiante al interactuar con esos recursos.

La conexión entre ambos planos se materializa mediante tres génesis:

| Génesis | Función en este estudio | Indicadores previstos |
|---|---|---|
| **Semiótica** | Coordinar registros algebraico, gráfico y numérico (A–G–N) en torno al producto de matrices. | Traducción matriz ↔ transformación de figuras, uso de tablas de coordenadas. |
| **Instrumental** | Transformar GeoGebra de simple artefacto a **instrumento cognitivo**. | Creación y uso autónomo de deslizadores, comando *AplicaMatriz*, exploraciones previas al cálculo. |
| **Discursiva** | Construir y comunicar argumentos sobre el efecto geométrico de las matrices. | Explicaciones orales o escritas que justifiquen rotaciones, cizalladuras o reflexiones. |

El análisis de la génesis instrumental se apoya en la Aproximación Instrumental (Rabardel, 1995), que diferencia entre artefacto (herramienta disponible) e instrumento (herramienta interiorizada por el estudiante mediante esquemas de uso). Según Drijvers (2013), esta transformación requiere procesos de instrumentalización (uso adaptado del software) e instrumentación (modificación del pensamiento matemático a partir del uso). Para el producto de matrices, GeoGebra actúa como artefacto clave al visualizar transformaciones geométricas y explorar la estructura de composiciones lineales, potenciando la construcción conceptual.

Complementariamente, se considera la Teoría de los Registros de Representación Semiótica de Duval (2006), que plantea que el aprendizaje significativo en matemáticas exige la conversión entre distintos registros (algebraico, gráfico, numérico y verbal). En el contexto del producto de matrices, se requiere que el estudiante relacione expresiones algebraicas con transformaciones gráficas y tablas de coordenadas, lo que constituye un eje central en las tareas propuestas.

**Marco metodológico**

Metodológicamente, este estudio adopta un diseño cualitativo de estudio de caso único (Korstjens & Moser, 2017) con 10 estudiantes de quinto de secundaria. Este enfoque es idóneo para analizar la evolución diacrónica del ETM mediante observaciones sistemáticas y análisis de producciones estudiantiles (tareas, registros GeoGebra), asegurando rigor a través de la triangulación de datos. Así se cumplen criterios de credibilidad, al tiempo que se genera transferibilidad analítica mediante la descripción específica de cómo la mediación tecnológica modula el ETM durante el aprendizaje de productos matriciales, revelando la transformación de sus componentes en interacciones prácticas.

La intervención didáctica se estructuró en dos fases secuenciales (preparación e implementación del ETM), distribuidas en cuatro sesiones:

Fase 1: Preparación del ETM

- Sesión 0 (Inducción instrumental): Familiarización operativa con GeoGebra para garantizar dominio básico del artefacto.
- Sesión 1 (Construcción algebraica): Deducción del algoritmo de multiplicación matricial y análisis de propiedades, priorizando la génesis semiótica en el plano simbólico.

Fase 2: Implementación del ETM mediado por tecnología

- Sesión 2 (Activación geométrica manual): Taller de transformaciones lineales mediante multiplicación matriz-vector, graficadas en papel milimetrado para consolidar la conexión álgebra-geometría sin mediación digital.
- Sesión 3 (Instrumentación tecnológica): Laboratorio dinámico con GeoGebra usando deslizadores y el comando AplicaMatriz, diseñado para promover la génesis instrumental y discursiva mediante exploración guiada de composiciones lineales."

El proceso de análisis se estructuró en tres etapas secuenciales para caracterizar la evolución del Espacio de Trabajo Matemático (ETM) personal durante la transición del producto de matrices de operación mecánica a objeto geométrico instrumentado:

- Primero, se aplicó una rúbrica rápida individual evaluando cuatro criterios por matriz (P-escala, Q-rotación, R-reflexión): 1) Corrección del cálculo fila×columna (A), 2) Coherencia numérica de coordenadas (N), 3) Precisión gráfica (G), y 4) Explicación del efecto geométrico (D). Esto generó una matriz de doble entrada que cuantificó los aciertos por habilidad.
- Segundo, se calcularon porcentajes de éxito por criterio y matriz, sintetizados en una tabla compacta de 12 indicadores. Esta permitió visualizar patrones de coordinación entre registros algebraicos, numéricos, gráficos y discursivos (A-N-G-D).

- Tercero, se desarrollaron viñetas micro-ilustrativas mediante contraste de casos: un ejemplo exitoso y otro problemático (respaldados por evidencia visual), descritos en 3 a 4 líneas para evidenciar tanto la circulación fluida entre registros como los bloqueos en la génesis instrumental.

Este abordaje tripartito reveló cómo la mediación tecnológica transforma la comprensión del producto matricial, documentando su reconstrucción como objeto explorado, argumentado e internalizado en el ETM estudiantil.

## Desarrollo y Resultados

La intervención se desplegó en dos fases secuenciales y cuatro sesiones, combinando trabajo manual y exploración digital para favorecer la transición desde un tratamiento puramente algorítmico del producto matricial hacia un Espacio de Trabajo Matemático (ETM) integrado. La Tabla 1 resume la secuencia didáctica, génesis predominantes y evidencias clave.

**Fase de preparación: Fundamentos del ETM**

- Sesión 0 (Inducción instrumental)

    8 de los 10 estudiantes construyeron matrices vinculadas a deslizadores en GeoGebra sin asistencia, evidenciando una génesis instrumental incipiente. Este dominio operativo básico fue crucial para reducir la carga técnica en sesiones posteriores.

- Sesión 1 (Construcción algebraica)

    El 70% de los estudiantes resolvió íntegramente productos matriciales en papel mediante el algoritmo fila×columna. Sin embargo, el trabajo se limitó al registro algebraico (A), sin articulación con otros planos del ETM, confirmando un tratamiento fragmentario inicial.

**Fase de implementación: Integración tecnológica del ETM**

- Sesión 2 (Taller manual de transformaciones)

    En papel milimetrado, los estudiantes aplicaron matrices de escalamiento $P = \begin{pmatrix} 2 & 0 \\ 0 & 2 \end{pmatrix}$, rotación $Q = \begin{pmatrix} 0 & -1 \\ 1 & 0 \end{pmatrix}$ y reflexión $R = \begin{pmatrix} 1 & 0 \\ 0 & -1 \end{pmatrix}$. Destacaron:

    **Articulación A-N-G**: El estudiante E-03 argumentó: *"El 2 en la diagonal duplica las dimensiones preservando la forma"*, coordinando registros algebraicos (A), numérico (N) y gráfico (G).

    **Dificultades conceptuales**: Solo 2 de 10 estudiantes intuyeron modificaciones para giro horario, y ningún estudiante dedujo cómo generar reflexión horizontal.

    **Hallazgo***:* La actividad desencadenó génesis semiótica activa y discursiva emergente mediante justificaciones basadas en efectos geométricos.

- Sesión 3 (Laboratorio dinámico con GeoGebra)
  Usando el comando AplicaMatriz y deslizadores $(a, b, c, d)$:

  **Génesis instrumental avanzada**: 70% manipuló parámetros independientemente, internalizando GeoGebra como instrumento cognitivo.

  **Coordinación consolidada**: Relacionaron variaciones numéricas con efectos geométricos en tiempo real (ejemplo *"Al mover $a$, la figura se estira en el eje x"*).

  **Génesis discursiva sólida**: 7 de 10 estudiantes formularon reglas precisas como: "*b causa inclinación izquierda-derecha; c controla inclinación arriba-abajo*".

Los resultados revelan una progresión sostenida en las tres génesis del ETM: desde un enfoque fragmentado (Sesión 1: solo álgebra) hasta una comprensión integrada de transformaciones lineales (Sesión 3: articulación A-N-G mediada por tecnología). Este avance confirma que el trabajo manual preparó una instrumentalización significativa de GeoGebra, permitiendo:

- La consolidación de la génesis semiótica mediante coordinación de registros,
- La internalización tecnológica (génesis instrumental), y
- La construcción de explicaciones matemáticas robustas (génesis discursiva).

**Tabla 1**

*Resumen de la secuencia didáctica y evolución de las génesis del ETM*

| Fase | Sesión | Propósito didáctico | Génesis predominante | Evidencias clave |
| --- | --- | --- | --- | --- |
| Preparación | 0. Inducción a GeoGebra | Reducir carga técnica | Instrumental incipiente | 8/10 vinculan deslizadores sin ayuda |
|  | 1. Construcción algebraica | Automatizar algoritmo fila×columna | Semiótica inicial (solo registro A) | 70% resuelve productos matriciales |
| Implementación | 2. Taller manual | Coordinar registros A→N→G | Semiótica activa + Discursiva emergente | E-03: *"El 2 duplica dimensiones"* |
|  | 3. Laboratorio dinámico | Usar GeoGebra como instrumento | Todas consolidadas: Instrumental + Semiótica + Discursiva | 7 de 10 verbalizan reglas de transformación |

## Análisis o discusión de resultados

Los hallazgos revelan cómo la secuencia didáctica basada en el Espacio de Trabajo Matemático (ETM) facilitó una activación progresiva y articulada de las tres génesis (Kuzniak, 2011), transformando la comprensión del producto matricial. Este análisis se estructura en tres dimensiones clave:

1. Génesis semiótica: Articulación multimodal de registros

   La coordinación algebraico-numérico-gráfica (A-N-G) en la Sesión 2 marcó un punto de inflexión: al calcular productos matriciales y representar gráficamente transformaciones en papel milimetrado, los estudiantes trascendieron el enfoque algorítmico tradicional criticado por Hidayanti (2020). Casos como el de E-02, quien vinculó el escalar 3 en $\begin{pmatrix} 3 & 0 \\ 0 & 3 \end{pmatrix}$ con el triplicado de dimensiones evidencian cómo el trabajo manual convirtió matrices en objetos de significado visual, sentando las bases para la modelación geométrica.

2. Génesis instrumental: Internalización tecnológica como herramienta cognitiva

   La Sesión 3 demostró la transición de GeoGebra de artefacto a instrumento (Rabardel, 1995): mediante AplicaMatriz y deslizadores, los estudiantes no solo visualizaron efectos, sino que desarrollaron esquemas de uso reflexivo (Drijvers, 2013). Interpretaciones como "b modifica la inclinación en x" revelan una apropiación funcional de parámetros matriciales $(a, b, c, d)$, donde la tecnología operó como extensión cognitiva para explorar relaciones estructura-efecto.

3. Génesis discursiva: Del lenguaje intuitivo a la argumentación formal

   La emergencia de verbalizaciones precisas, desde inferencias sobre rotaciones (E-01) hasta descripciones paramétricas en Sesión 3 señala la consolidación de la génesis discursiva. Este tránsito desde intuiciones locales ("cambiar signos altera el giro") hacia reglas generales ("c controla inclinación vertical") manifiesta la construcción de un lenguaje matemático funcional, indicador clave de internalización conceptual en el ETM.

**Síntesis teórica: Integración dialéctica de génesis**

Los resultados validan el principio de Kuzniak (2011), solo el desarrollo equilibrado de las tres génesis genera comprensión profunda. La secuencia logró esta sinergia mediante:

- Fase manual (Sesión 2): Activó semiótica y discurso emergente sin mediación tecnológica,
- Fase digital (Sesión 3): Potenció instrumentalización y generalización discursiva mediante exploración guiada.

Esta dialéctica corrobora que la tecnología no sustituye el razonamiento, sino que lo amplifica cuando se inserta en diseños pedagógicos intencionados (Drijvers, 2013).

**Implicaciones pedagógicas y alcance**

Más allá del producto matricial, la secuencia desarrolló competencias transversales:

- Pensamiento relacional: Conexiones A-N-G sustentadas en visualización,
- Alfabetización tecnocrítica: Uso reflexivo de GeoGebra para verificar conjeturas,
- Comunicación matemática: Argumentación basada en evidencias observables.

Estos resultados proponen un modelo replicable para enseñar álgebra lineal: integrar manipulación concreta, exploración tecnológica y construcción discursiva, donde cada fase consolida dimensiones complementarias del ETM.

## Conclusiones y recomendaciones

Esta investigación demuestra que la articulación coherente de los registros algebraico, gráfico y numérico en la enseñanza del producto de matrices es viable cuando se implementan secuencias didácticas que exigen explícitamente dicha integración. El diseño híbrido (fase manual + fase digital) resultó fundamental para que los estudiantes trascendieran una comprensión operativa del algoritmo hacia una conceptualización geométrica de las transformaciones lineales, evidenciando la evolución del Espacio de Trabajo Matemático (ETM) mediante tres hallazgos centrales:

1. El trabajo manual como base diagnóstica: Las tareas con papel y lápiz (Sesión 2) no solo desarrollaron intuiciones geométricas autónomas, sino que identificaron obstáculos semióticos críticos (ejemplo: dificultad para generalizar reflexiones horizontales), proporcionando una base cognitiva para la posterior mediación tecnológica.
2. GeoGebra como catalizador de procesos superiores: La instrumentalización del software (Sesión 3) potenció la visualización dinámica, facilitó la generación de conjeturas sobre parámetros matriciales (a, b, c, d) y consolidó la argumentación matemática mediante descripciones funcionales (ejemplo: "b modifica la inclinación en el eje x").
3. El lenguaje como indicador de apropiación conceptual: Las verbalizaciones estudiantiles evidenciaron una progresiva internalización del producto matricial como objeto de exploración (no solo cálculo), manifestando la consolidación de la génesis discursiva del ETM.

Si bien los resultados son promisorios, se reconocen dos limitaciones estructurales:

- El tamaño muestral reducido (10 estudiantes) y el contexto educativo específico restringen la transferibilidad estadística de los hallazgos.
- La brevedad de la intervención (4 sesiones) imposibilita evaluar la permanencia de los aprendizajes a mediano plazo.

Para superar estas barreras y profundizar los alcances, se proponen ampliaciones en el estudio:

- Incorporar la exploración explícita de la no conmutatividad (ejemplo: comparar efectos de $P.Q$ vs $Q.P$.
- Analizar el determinante como indicador de cambios de escala, orientación e invertibilidad.

En conclusión, este trabajo reafirma que la enseñanza del álgebra lineal debe fundamentarse en diseños que articulen representaciones, herramientas y discurso matemático. La tecnología, lejos de reemplazar el razonamiento lo potencia cuando se subordina a objetivos pedagógicos claros, se combina con fases de manipulación no digital, y exige coordinación activa entre registros. La evolución observada en las génesis del ETM corrobora que solo mediante esta integración multidimensional el producto de matrices deja de ser un procedimiento abstracto para convertirse en un concepto explorado, instrumentado y argumentado.